\documentclass[12pt]{amsart}
\def\thZ{ \mathcal{H}}

\def\JJ{J_\alpha}

 \def\NN{\mathcal{N}}

\def\JJ{J_\alpha}

\usepackage[square,compress,comma, numbers,sort]{natbib}
\usepackage[colorlinks=true, citecolor=red, linkcolor=blue]{hyperref}
\usepackage{amsfonts,mathtools}

\allowdisplaybreaks[4]

\def\tTh{\widetilde{\Theta}}

\def\tX{\widetilde{X}}
\def\tZ{\widetilde{Z}}

\usepackage{amssymb}
\usepackage{color}
\newcommand{\abs}[1]{\left\lvert #1 \right\rvert}
\newcommand{\Abs}[1]{ \biggl \lvert #1 \biggr \rvert}

\DeclarePairedDelimiterXPP\pk[1]{\mathbb{P}}\{ \}{}{ #1}
\DeclarePairedDelimiterXPP\E[1]{\mathbb{E}}\{ \}{}{	#1}

\def\FRE{\mbox{Fr\'{e}chet }}

\def\toprob{\overset{p}\rightarrow}
 \usepackage{xparse}% http://ctan.org/pkg/xparse

\NewDocumentCommand{\ceil}{s O{} m}{%
  \IfBooleanTF{#1} % starred
    {\left\lceil#3\right\rceil} % \ceil*[..]{..}
    {#2\lceil#3#2\rceil} % \ceil[..]{..}
}
\NewDocumentCommand{\floor}{s O{} m}{%
  \IfBooleanTF{#1} % starred
    {\left\lfloor#3\right\rfloor}
    {#2\lfloor#3#2\rfloor}
}

\definecolor{c20}{rgb}{0.,0.7,0.}
\definecolor{c30}{rgb}{0.,0.,1.}
\definecolor{c40}{rgb}{1,0.1,0.7}
\definecolor{c50}{rgb}{1,0,0}
\definecolor{c60}{rgb}{1,0.9,0.1}
\definecolor{c70}{rgb}{0.50,1.00,0.00}

\def\bE#1{{\textcolor{c30}{#1}}}

\def\bE#1{#1}

\def\Ecc#1{{\textcolor{c50}{#1}}}
\def\Ecc#1{#1}

\def\kd#1{{\textcolor{black}{#1}}}
\def\kdd#1{{\textcolor{black}{#1}}}

\numberwithin{equation}{section}
\newtheorem{theo}{Theorem}[section]
\newtheorem{sat}[theo]{Proposition}
\newtheorem{de}[theo]{Definition}
\newtheorem{lem}[theo]{Lemma}

\newtheorem{example}[theo]{Example}
\newtheorem{korr}[theo]{Corollary}
\newtheorem{remark}[theo]{Remark}

\numberwithin{equation}{section}

\newcommand{\prooftheo}[1]{ \textsc{Proof of Theorem} \ref{#1} }

\newcommand{\prooflem}[1]{\textsc{Proof of Lemma} \ref{#1}}

\newcommand{\QED}{\hfill $\Box$}

\newcommand{\COM}[1]{}

\def\IF{\infty}

\newcommand{\R}{\mathbb{R}}
\newcommand{\inr}{\in \R}

\newcommand{\BQN}{\begin{eqnarray}}
\newcommand{\EQN}{\end{eqnarray}}
\newcommand{\BQNY}{\begin{eqnarray*}}
\newcommand{\EQNY}{\end{eqnarray*}}
\def\ldot{, \ldots,}

\newcommand{\limit}[1]{\lim_{#1 \to   \infty}}

\def\todis{\overset{d}\rightarrow}

\def\kk#1{\textcolor{cyan}{#1}}
\def\kk#1{{#1}}

\def\bqny#1 {\begin{eqnarray*} #1 \end{eqnarray*}}
\def\bqn#1{\begin{eqnarray} #1 \end{eqnarray}}

\newcommand{\BS}{\begin{sat}}
\newcommand{\ES}{\end{sat}}
\newcommand{\BT}{\begin{theo}}
\newcommand{\ET}{\end{theo}}
\newcommand{\BK}{\begin{korr}}
\newcommand{\EK}{\end{korr}}

\newcommand{\BEX}{\begin{example}}
\newcommand{\EEX}{\end{example}}

\newcommand{\BD}{\begin{de}}
\newcommand{\ED}{\end{de}}

\newcommand{\BIT}{\begin{itemize}}
\newcommand{\EIT}{\end{itemize}}
\newcommand{\BDI}{\begin{description}}
\newcommand{\EDI}{\end{description}}

\newcommand{\BRM}{\begin{remark}}
\newcommand{\ERM}{\end{remark}}

\newcommand{\BEL}{\begin{lem}}
\newcommand{\EEL}{\end{lem}}

\newcommand{\nelem}[1]{{Lemma \ref{#1}}}
\newcommand{\neprop}[1]{{Proposition \ref{#1}}}
\newcommand{\netheo}[1]{{Theorem \ref{#1}}}

\def\TT{\mathcal{T}}
\def\TT{\R }

\def\Z{\mathbb{Z}}

\def\JJ{\mathcal{S}(Z)}
\def\JJd{\mathcal{S}(Z)}

\def\TT{\R}
\def\intT{\int_{\TT}}

\begin{document}

\title[Supremum of Max-Stable Processes \& Pickands Constants]
{}

\author{Krzysztof D\c{e}bicki}
\address{Krzysztof D\c{e}bicki, Mathematical Institute, University of Wroc\l aw, pl. Grunwaldzki 2/4, 50-384 Wroc\l aw, Poland}
\email{Krzysztof.Debicki@math.uni.wroc.pl}

\author{Enkelejd  Hashorva}
\address{Enkelejd Hashorva, Department of Actuarial Science, University of Lausanne,\\
Chamberonne 1015 Lausanne, Switzerland}
\email{Enkelejd.Hashorva@unil.ch}

\bigskip

\date{\today}
 \maketitle

\begin{quote}

{\bf Abstract:} Let $X(t),t\in \R$ be a stochastically continuous stationary max-stable process with
\FRE marginals $\Phi_\alpha, \alpha>0$ and set
 $M_X(T)=\sup_{t \in   [0,T]} X(t),T>0$.  In the light of the seminal articles  \cite{Genna04,Genna04c},  it follows that $A_T=M_X(T)/T^{1/\alpha}$ converges in  distribution as $T\to \IF$ to $\thZ^{1/\alpha} X(1)$, where $\thZ$ is the Pickands constant corresponding to the spectral process $Z$ of $X$. In this contribution we derive explicit formulas for $\thZ$ in terms of $Z$ and show necessary and sufficient conditions for its positivity. From our analysis it follows that $A_T^\beta,T>0$ is uniformly integrable for any $\beta \in (0,\alpha)$.
 For Brown-Resnick $X$ we show the validity of the celebrated Slepian inequality
 and discuss the finiteness of Piterbarg constants.
\end{quote}

{\bf Key Words:}
Max-stable process; spectral tail process;  Gaussian processes with stationary increments; L\'evy processes;  Pickands constants; Piterbarg constants;  Slepian inequality; growth of supremum.

{\bf AMS Classification:} Primary 60G15; secondary 60G70\\

\section{Introduction}
Let $X(t), t\in \TT$ be a stochastically continuous stationary max-stable process with \FRE marginals $\Phi_\alpha(x)=e^{-x^{-\alpha}}, x>0,\alpha>0$.
Here max-stable means that the finite dimensional distributions (fidi's) of $X$ are max-stable multivariate distributions, see e.g., \cite{deHaan, kab2009}. \bE{By a key theorem of  de Haan \cite{deHaan}, it is well-known
that $X$ can be represented (in distribution)  as}
\bqn{\label{eq1}
	X(t) =  \max_{i\ge 1} P_i Z^{(i)}(t), \quad t\in \TT,
}
where $\Pi= \sum_{i=1}^\IF \varepsilon_{P_i}$ is a Poisson point process (PPP) on $[0,\IF) $ with intensity
$\alpha x^{-\alpha- 1} dx$ independent of $Z^{(i)}$'s which are independent copies of a random process $Z(t),t\inr$ which shall be referred to as the spectral process. The assumption that $X(t)$ has distribution $\Phi_\alpha$ implies that
%$(we follow )the following It is well-known In the following $Z(t),t\in \TT$ is a non-negative process satisfying
$ \E{Z^\alpha(t)}=1, \forall t\in \TT$
%and $X$ has spectral process $Z$, i.e., we  have  the de Haan representation, see e.g., \cite{deHaan,dom2016} (below $\EQD $ means equality in law)
  Since we consider here only stationary max-stable processes, adapting the terminology of \cite{kab2009}
  we shall call the spectral process $Z$ a  {\it Brown-Resnick} stationary process.

The assumption that $X$ is stochastically continuous implies that it has a separable and measurable version, see e.g., \cite{Doob};  the same holds for $Z$, see \cite{dom2016}. Therefore in the following we suppose that both $X$ and $Z$ are jointly measurable and separable. Hereafter  we shall assume further that $X$ has locally bounded sample paths, and thus by \eqref{eq1} $Z$ also has locally bounded  sample paths.
According to  \cite{dom2016},  this assumption is important for conditions that guarantee the existence of a dissipative  Rosi\'nski (or also called a mixed moving maxima)  representation of $X$.

By separability and the \bE{local} boundedness of the sample paths of $X$ and $Z$ both
$M_X(T)=  \sup_{t \in   [0,T]} X(t)$ and $M_{Z^\alpha}(T)$ are well-defined and finite random variables for any $T>0$.
By \eqref{eq1}, given
$t_i\inr, x_i \in (0,\IF), i\le k$ we have (see e.g., \cite{DM})
\bqn{\label{foll}
	- \ln \pk*{\kd{\forall_{1\le i \le k }} X(t_i) \le  x_i}=
  \E[\Big]{ \max_{1 \le i \le \kd{k}} Z^\alpha (t_i)/x_i^\alpha }.
}
 By the measurability of $Z$, for any $T>0$ using Fubini theorem
$$\E*{\int_0^T Z^\alpha(t)\lambda(dt)}=\int_0^T \E{Z^\alpha(t)} \lambda(dt)= T,$$
\kd{with $\lambda(\cdot)$ the Lebesgue measure on $\R$}. Further the separability assumption and the above considerations imply 
% we have % then \bE{using further \eqref{nat}} we have
$$ \pk{ M_X(T) \bE{<} \IF}= \pk{ M_{Z^\alpha}(T) \bE{<} \IF}=1$$
and
\bqn{ - \ln \pk*{ M_X(T) \le (Tx)^{1/\alpha} }
	=   \frac{\E*{ M_{Z^\alpha}(T)}}{Tx}, \quad \forall T,x>0.\label{fred}
}
Hence  $\E{ M_{Z^\alpha}(T)},T>0$ does not depend on the particular choice of the \kdd{spectral process} $Z$ but only on $X$. By the stationarity of $X$ it follows that
$$\E*{ M_{Z^\alpha}(T)}= \E*{ \sup_{ S \le t \le S+T } Z^\alpha(t)} \bE{ \in (  0, \IF)}$$
 for any  \kd{$S\inr,T>0$}. Hence $\E{ M_{Z^\alpha}(T)}, T>0$ is  sub-additive and consequently,  Fekete lemma yields
\bqn{\label{spilman}
	\thZ\kd{:=} \limit{T}\frac{1}{T}\E[\big]{ M_{Z^\alpha}(T)}=
\inf_{T> 0} \frac{1}{T}\E[\big]{ M_{Z^\alpha}(T)}  \le  \E[\big]{ M_{Z^\alpha}(\kd{1})} \bE{\in (0, \IF)}.
}
Moreover, from the above we conculde that $\thZ$ does not depend on the particular choice of the spectral tail process $Z$ but only on the stationary max-stable process $X$. Referring to \cite{SBK}, $\thZ$ is the so-called {\it generalised Pickands} constant defined with respect to \bE{some} Brown-Resnick stationary process $Z$. In \cite{debicki2002ruin} $\thZ$ is introduced for the log-normal process
\bqn{\label{GaussStat}
	Z(t)&=& e^{ B(t)- \sigma^2(t)/2}, \quad t\inr,
}
where $B(t),t\inr $ is a centered Gaussian process with stationary increments, continuous sample paths  and variance function $\sigma^2$ \bE{which does not vanish in any interval of $\R$}.\\
Taking $B$ to be a fractional Brownian motion (fBm) with self-similarity Hurst index $\alpha/2\in (0,1]$,
\kd{we get that} $\thZ$ is the {\it classical Pickands} constant, see e.g., \cite{PicandsA, Pit96, SBK}. The only known values of $\thZ$  are 1 and $1/\sqrt{\pi}$ corresponding to $\alpha=1$ and $\alpha=2$, respectively.
 In the case of L\'evy processes the Pickands constant $\thZ$ appears explicitly in many contributions, see e.g.,
 \cite{ KW,SBK} and references therein. Moreover, for the discrete-time case $X(t),t\in \Z$ we have that $\thZ$
 (introduced similarly as for the continuous-time, see \cite{KDEH1}) is  the extremal index of the stationary
 time series $X(t),t\in \Z$. In that context, it has been also studied in \cite{MR2453345} using the spectral
 representation of $X$. A considerable amount of research is dedicated to  calculation and estimation of the extremal
  index of regularly varying time series, see e.g., \cite{Hrovje} and the reference therein. \\
The main question that arises for Pickands constants $\thZ$ is:\\
   {\bf Q1}: Under what conditions are these constants positive or equal to 0?\\
For \bE{a  } stationary  max-stable process $X$ this question is  partially answered in \cite{SBK} \bE{when}  $Z$ is such that $Z(0)=1$ almost surely. The case that  $Z(0)$ is  a non-negative random variable  in treated in \cite{Htilt}. Specifically, the positivity of $\thZ$ has been shown under the assumption that
\bqn{
	\label{eqZ} \pk*{  \JJ < \IF}=1, \quad \JJ:=\int_\R Z^\alpha(t) \lambda(dt).
}
In view of \cite{dom2016}   since $X$ has locally bounded sample paths, then under \eqref{eqZ} $X$ has a dissipative Rosi\'nski  representation which is equivalent with $X$ being generated by a non-singular dissipative  flow, see
\cite{kab2009a,StoevWangSPL,Roy, dom2016} for more details.
As shown in  \cite{SBK, Htilt}, if the spectral process $Z$ has c\`adl\`ag sample paths and \eqref{eqZ} holds, then
\bqn{ \label{gjge}\
	\thZ \ge  \E*{\frac{ \sup_{t \in  \R} Z^\alpha (t)}{\JJ}}
,}
with equality shown under some technical  assumptions for both
Gaussian and L\'evy spectral processes $Z$.  The investigation therein was motivated by \cite{DiekerY, DM}. The former contribution showed that
\kdd{\eqref{gjge} holds with equality for $B$ in \eqref{GaussStat} being an fBm.}
Since  $\thZ$ in \eqref{spilman} is  defined as a limit, it turns out that the explicit calculation of $\thZ$ is for general $Z$ too difficult. However, if  \eqref{gjge} holds with equality, then $\thZ$ being an expectation,
can be efficiently simulated, see e.g.,  \cite{DiekerY}.\\
  An interesting question that arises here is: \\
{\bf Q2}: Does \eqref{gjge} hold with equality for general Brown-Resnick stationary $Z$?\\
Clearly, if $\thZ=0$, then \eqref{spilman}  means the convergence in probability
\bqn{ A_T:=\frac{M_X(T)}{T^{1/\alpha}} \toprob 0, \quad T\to \IF,
\label{at1}}
whereas when $\thZ>0$ we have the convergence in distribution
\bqn{\label{at2} A_T \todis \thZ^{1/\alpha} X(1), \quad T\to \IF.}
For $X$ being a symmetric $\alpha$-stable ($S\alpha S$)  stationary process  with $\alpha \in (0,2)$ the above convergence has been shown in the seminal articles \cite{Genna04,Genna04c}, see the recent contributions  \cite{MR2384479, Roy0,Roy1, Roy2} for related results and new developments.\\
The findings of \cite{Genna04,Genna04c} are important for the max-stable processes too,
which is already pointed out in \cite{MR2453345} for discrete max-stable processes.
Indeed, using the link between max-stable and $S\alpha S$ processes established
in \cite{kab2009a} and \cite{WangStoev} independently, it follows that
when $X$ is generated by a non-singular conservative flow, which by \cite{dom2016} (under the assumption of locally boundedness of sample paths of $X$) is equivalent with  $\pk*{ \JJ = \IF}=1,$   then we have
\bqn{\thZ=0. \label{wessen}
}
Note that \eqref{wessen} holds also when we consider the discrete case $X(t),t\in \Z$, which can be shown  for instance by utilising the expression of Pickands constant (which in this case coincides with the extremal index, \cite{KDEH1,Htilt}) derived in \cite{MR2453345}. \\
We conclude that $\thZ$ is  positive if and only if
$$\pk*{ \JJ = \IF}<1.$$
Hence according to our argumentation above question {\bf Q1} has a simple answer.
Namely, if $X$ is a stationary max-stable process with locally bounded sample paths,
then (by \cite{dom2016}),
$\thZ>0$ if and only if
\bqn{\label{bask}
	\pk{ \JJ< \IF}>0.
	}
Clearly, the convergence in probability in \eqref{at1} implies that
$A_T^\beta \toprob 0$ as $T\to \IF$ for any $\beta \in (0,\IF)$
and a \kd{similar} implication holds for $A_T^\beta$ when  \eqref{at2} is satisfied. For $X$ being an  $S\alpha S$ random field the recent contribution \cite{RoyY} strengthened those convergences to that  of $\E{A_T^\beta}$ for $\beta \in (0,\alpha)$, i.e., showing the uniform integrability of $A_T^\beta$ whenever $\beta \in (0,\alpha)$. The case that $X$ is a stationary max-stable random process  is easier to deal with, see \neprop{ProK} in Section 4. \\
Our main interest \kd{in this contribution} is the derivation of expressions for $\thZ$
in terms of the spectral process $Z$ that appears in the de Haan
representation \eqref{eq1}. In particular, motivated by {\bf Q2}, we
show that \eqref{gjge} (or a modification of it) holds with equality  under
\eqref{bask} without further assumptions.  Recall that so far it is only known that the inequality in \eqref{gjge} holds for $X$ having a dissipative Rosi\'nski representation. \\
 As already shown in \cite{SBK,Htilt}, different representations for $\thZ$ relate to different dissipative Rosi\'nski  representations of $X$. Therefore, our analysis is also concerned with such representations for $X$.

Our study  of Pickands constants (together with the criteria for its positivity) allows us to
investigate the growth of the expectations of $M_X(T)$  and $M_{Z^{\alpha}}(T)$ as $T\to \IF$.
The latter can be investigated under the further assumption of the
Brown-Resnick model, i.e., \kdd{when} $Z$ is a log-normal process.  Moreover, for the Brown-Resnick model
an extension of the celebrated Slepian inequality is possible, see \netheo{th.pick} below.

 Organisation of the paper: Our main results are displayed in Section 2 followed by discussions and some extensions \bE{presented} in Section 3. Proofs are postponed to Section 4; an Appendix concludes this contribution.

\def\AA{\mathcal{A}}

 \section{Main Results}
Let $X,Z$ be \kk{(as in the Introduction)}  jointly measurable, separable and with
 locally bounded sample paths.  By the measurability of $Z$ we have that $\JJ=\int_\R Z^\alpha(t) \lambda(dt)$ is a random variable in $\R \cup \{+ \IF\}$, see \cite{Doob}.  Write next $\E{A; B}$ instead of $\E{A \mathbb{I}(B ) }$ for an event $B$ with $\pk{B}>0$. Fixing $T>0$ we have the following splitting formula
 \bqn{\label{split} \E{M_{Z^\alpha}(T)} =  \E{M_{Z^\alpha}(T); \JJ < \IF}
 +  \E{M_{Z^\alpha}(T); \JJ= \IF}.
}
If $\pk{ \JJ=\IF}>0$,  \kd{then} by \cite{DombKabB}[Lem 16] the random process $Z_D$ defined by
$$ Z_D(t):=Z(t) \mathbb{I}( \JJ= \IF),\quad t\inr$$
 is Brown-Resnick stationary. Since $Z_D$ has also locally bounded sample paths and $\mathcal{S}(Z_D)= \IF$ almost surely,  the corresponding max-stable process $X_D$ is generated by a non-singular conservative flow. Moreover, under  \eqref{bask}
$$Z_C(t)=
Z(t) \mathbb{I}( \JJ< \IF)$$
 is also a Brown-Resnick stationary process which is generated by a
 non-singular dissipative flow. In order to omit technical details, we refer the reader to
 \cite{Roy, dom2016} for details on conservative and dissipative parts of max-stable processes.
 Consequently, by the discussions in the Introduction, condition \eqref{bask} implies that
%if $X$ is not generated by a non-singular  conservative flow
\bqn{\label{zetC} \thZ=\mathcal{H}_{Z_C}> 0,
}
where $\mathcal{H}_{Z_C}$ is the Pickands constant with respect to spectral process $Z_C$.
{We remark that \eqref{zetC} is new and not available in the literature so far.}\\
 In view of \eqref{zetC},  in the following we can reduce our analysis by considering only the case that $X$
 is generated by a non-singular dissipative flow. In view of \cite{dom2016} this is equivalent with $X$ having  a dissipative Rosi\'nski  representation i.e., for some  non-negative random process (called also random shape function) $L(t),t\inr$ which is continuous in probability (we can consider here therefore $L$ to be jointly measurable and separable)  and for some $c>0$ we have the representation (in distribution)
\bqn{ \label{ali}
X(t) =  \max_{i\geq 1}  P_i L^{(i)}(t - T_i) , \quad t\in  \TT,
}
where $\sum_{i=1}^\infty \epsilon_{(P_i, T_i)}$ is a PPP on $[0,\IF) \times \TT$ with intensity $ c\cdot \lambda(dt) \cdot  \alpha x^{-\alpha -1} dx$, independent of $L^{(i)}$'s which are independent copies of $L$.\\
 By \eqref{ali}  for any random variable $\NN$ with density $p(t)>0, t\in \TT$
\bqn{\label{piratenQ}
	Z =  (c / p(\NN)) ^{1/\alpha}B^\NN L
}
is a valid spectral process for $X$, where  $\NN$ is independent of $L$ and  $B^t L(\cdot)= L(\cdot- t),t\inr$.

\BT \label{proviant}Let $X(t),t\in \R$ be a stochastically continuous max-stable \bE{stationary} process with de Haan representation \eqref{eq1}. % \bE{and spectral process $Z$ satisfying \eqref{nat}.}
 If $X$ has locally bounded sample paths  and condition  \eqref{eqZ} holds,
then there exists  some jointly measurable and separable non-negative random shape function $L$  such that \eqref{piratenQ} holds and moreover
 we have
\bqn{ \label{zaub}
\thZ = \frac{ \E{  \sup_{t \in \R} L^\alpha (t)}}{ \E{\mathcal{S}(L)}}
 \in (0,\IF).
}
\ET

\BRM
\COM{
	i) If $X(t),t\in \delta \Z, \delta >0$ is a stationary max-stable process, then the Pickands constant is defined by (see e.g., \cite{SBK})
$$\limit{n} \frac 1 T \E*{
	\sup_{t \in   \delta \Z\cap [0,T] } Z^\alpha(i)}=\thZ^\delta.
$$  In view of \cite{dom2016} $X$ has a dissipative Rosi\'nski  representation if and only if
$\pk{\sum_{t \in \delta \Z} Z(t) < \IF}>0$.  Under this condition, as in the proof of \netheo{proviant}, it follows that the Pickands constant is given by
\bqn{
\thZ^\delta  = \frac{ \E{  \sup_{t \in \delta \Z } L^\alpha (t)}}{
	\E{ \delta \sum_{t \in \delta \Z} L^\alpha (t)}}  \in (0,\delta].
}
%Note that for the special case that $\pk{\sum_{t \in \delta \Z} Z(t) < \IF}=1$ the above formula follows from \cite{Htilt} or \cite{Hrovje}. \\
ii)
}
  When $X$ has a dissipative Rosi\'nski  representation, it is possible to construct $L$ such that
$\sup_{t \in   \R} L^\alpha (t)=1$ almost surely, see \cite{dom2016}. Hence by \eqref{zaub} for such random shape functions $L$ we have
\bqn{
\thZ = \frac{ 1}{ \E{\mathcal{S}(L)}}.
}
We shall show below that it is also possible to construct $L$ such that
$\mathcal{S}(L)=1$ almost surely, and thus by \eqref{zaub} we obtain \kd{an}
alternative formula, namely
\bqn{
\thZ = \E*{\sup_{t \in   \R} L^\alpha(t)}.
}
\ERM

For simplicity we shall assume in the following that both $X$ and $Z$ have c\`adl\`ag sample paths.
Let $D$ be the space of c\`adl\`ag functions $f: \bE{[0,\IF)}\to \R$ equipped with a metric $d$ which makes it  complete and separable, see e.g., \cite{MR838085} for details. Let $\mathcal{D}$ be the Borel $\sigma$-algebra on $D$ defined by this metric and let $\mu$ be a probability measure   given by (interpret $0:0$ as 1 below)
$$ \mu( A)= \E{ Z^{\alpha}(0) \mathbb{I}( Z/Z(0) \in A)}, \quad
A \in \mathcal{D}.$$
Since $(D,d)$ is a Polish metric space (complete and separable), we can determine a stochastic process $\Theta(t),t\in \R$ with c\`adl\`ag sample paths and probability law $\mu$;  refer to $\Theta$ as the spectral tail process.
 By \cite{Htilt}[Thm 4.1] all the fidi's of the  max-stable process $X$ are determined by $\Theta$. Namely, we have the following inf-argmax formula valid for $x_i$'s positive constants and $t_i$'s in $\R$
	\bqn{ %\label{argmax}
		- \ln \pk{ X(t_i)\le x_i, 1 \le i\le n}
%		&=& \E*{ \max_{1  \le i \le n}Z^\alpha(t_i)/x_i^\alpha}\\
	&=& 	\sum_{k=1}^n \frac 1 {x_k^\alpha } \pk[\Big]{\inf {\rm argmax}_{1 \le i \le  n} \Bigl( \frac{\Theta^\alpha (t_i-t_k)  }{x_i^\alpha }\Bigr) = k}. \label{argmax2}
}
Consequently, $\Theta$ defines $X$ and vice-versa, from $X$ we can calculate the fidi's of $\Theta$ by the generalised Pareto distributions of $X$, see \cite{Htilt}[Remark 6.4]. The next result gives an explicit construction for the random shape function $L$ and confirms  \eqref{gjge}.

\BT \label{th2} Under the setup of \netheo{proviant}, if further $X$ has c\`adl\`ag sample paths and \eqref{bask} holds, then
\bqn{ \thZ =
	\E*{ \frac{ \sup_{t \in  \R } \Theta^\alpha (t)}{\mathcal{S}(\Theta)}; \mathcal{S}(\Theta) \in (0,\IF)}
\in (0,\IF).
}
Moreover, if \eqref{eqZ} is valid, then
$X$ has dissipative Rosi\'nski  representation \eqref{piratenQ} with random shape function
\bqn{
	L(t)= \frac{\Theta(t)}{  (\mathcal{S}(\Theta))^{1/\alpha}},\quad  t\in \TT.
}
\ET

\BRM If $X$ is as in \netheo{th2} with c\`adl\`ag sample paths, then it follows \kk{straightforwardly} that
\kk{ $\thZ=0$} is equivalent with $ \pk{\mathcal{S}(Z) = \IF}= \pk{\mathcal{S}(\Theta) = \IF}=1$. \bE{Note further that since $\Theta(0)=1$ almost surely, then $ \mathcal{S}(\Theta)>0$ almost surely since  $\Theta$ has c\`adl\`ag sample paths.}
\ERM

{\bf Example 1}. Consider the Gaussian case with $Z$ as in \eqref{GaussStat},
where $B(t),t\inr $ is a centered Gaussian process with stationary increments, continuous sample paths  and variance function $\sigma^2$ \bE{that does not vanishes in compact intervals of $\R$.}

 We can assume without loss of generality (see \cite{kab2009}) that $\sigma(0)=0$. Hence $Z(0)=\bE{1}$  almost surely and for the corresponding spectral tail process $\Theta$ we simply have $\Theta= Z$. Since for this case $\alpha=1$, then under  \eqref{bask}
\bqny{ \thZ =
	\E*{ \frac{ \sup_{t \in  \R } Z   (t)}{\JJ}; \JJ \in (0,\IF)}
	\in (0,\IF).
}
In view of \cite{kab2009}, the following condition
\bqn{\label{caj}
	\liminf_{t\to \IF} \frac{\sigma^2(t)}{ \ln t}> 8
}
 implies \eqref{eq1}  and thus  $X$ has a dissipative Rosi\'nski  representation with random shape function $L(t)= Z(t)/\mathcal{S}(Z), t\inr$. Moreover
\bqn{\thZ =
\E*{ \frac{ \sup_{t \in  \R } Z (t)}{\JJ}} \in (0,\IF),
} which has been proved in  \cite{SBK}[Thm 2] under  \kk{some}
additional assumptions on the variance function \bE{of $B$.}

{\bf Example 2}. Stationary max-stable L\'evy--Brown--Resnick
processes $X$  have  spectral processes
$Z(t)= e^{W(t)},t\inr$   constructed from two independent L\'evy processes.
 Specifically, let  $\{B^+(t), t\geq 0\}$ be a L\'evy process with Laplace exponent $\Psi(\theta)=\ln \E{\exp{\bigl( }\theta B^+(1) {\bigr)}} $ being finite for $\theta =1$. Write $-W^{-}$ for another independent L\'evy process  with Laplace exponent
$$\ln \E{e^{ \theta W^{-}(1)}}= \Psi(1-\theta)-(1-\theta)\Psi(1).$$
Then we set $W(t)= W^{+}(t)\kd{:=B^+(t)-\Psi(1)t}, t\ge 0$, and $W(t)=W^{-}(-t)$ if $t< 0$.
In view of \cite{eng2014c} the max-stable process $X$ with unit \FRE marginals $\Phi_1(x)=e^{-1/x},x>0$ corresponding to the spectral process $Z$ is stationary. \bE{Bote that this fact is proved in a completely different context in \cite{Sasha}[Lem 1].} \\
By \cite{eng2014c} the L\'evy-Brown--Resnick process $X$ admits a dissipative Rosi\'nski  representation and thus \netheo{th2} and \cite{SBK}[Thm 3.2]   imply that (note that since $Z(0)=\bE{1}$ almost surely, then  $\Theta= Z$)
\bqn{
\thZ= \E*{ \frac{ \sup_{t \in  \R } Z(t)}{\mathcal{S}  (Z)}} =
  \frac{\underline{k}(0,1)}{\underline{k}'(0,0)} > 0,
}
where $\underline k$ is the bivariate Laplace exponent of the descending
ladder process corresponding to \kd{$W^+$}. \\
 If $B^+$ is a spectrally negative L\'evy process, we have the alternative formula
 % since the righ hand side above equals $\Psi'$  we have further %by  \cite{eng2014d}
$\thZ = \Psi'(1)>0$, which is already derived  in \cite{eng2014d}.
%,
%\kd{which agrees with the fact that}
%$$	\Psi'(1) = \frac{\underline{k}(0,1)}{\underline{k}'(0,0)} > 0.
%$$

\section{Discussions \& Extensions}
\subsection{Slepian inequality  for Brown-Resnick max-stable processes}
Slepian inequality is essential in the theory of extremes and sample path properties of Gaussian and related processes.
%Besides it is also useful in numerous fields of mathematics including  optimisation and number theory problems, see e.g., \cite{MR3059269}.\\
 A \kd{commonly} used version of Slepian inequality given for instance in \cite{GennaSlepian}[Thm 1.1] is as follows:\\
  If $B_1(t), B_2(t),t\inr $ are two centered Gaussian processes, then  for any $t_1 \ldot t_n \inr, n\ge 1$
  we have
\bqny{ \E*{\max_{1 \le i \le n} B_1(t_i) } \ge
	\E*{\max_{1 \le i \le n}  B_2(t_i) } ,
}
provided that for all $1 \le i \not=j \le n$
\bqn{\label{aq}
	\E*{ (B_1(t_i)- B_1(t_j))^2} \ge \E*{ (B_2(t_i)- B_2(t_j))^2}.
}	
\kdd{Moreover, in view of \cite{vitale2}[Eq. 6] (applied to $g(x)=e^x$)}
 for any real-valued function $f$
\bqn{ \label{tim}
		\E*{ \max_{1 \le i,j \le n}e^{{[B_1(t_i)- \kk{B_1}(t_j)]}- f(t_i)}} \ge
	\E*{ \max_{1 \le i,j \le n}e^{{[B_2(t_i)- B_2(t_j)]}- f(t_i)}}.
}	
%provided that $B(t_0)=0$ almost surely for some $t_0 \in \R$. \\
Let  $X_i(t),t\inr,i=1,2$ be max-stable processes with spectral processes  $Z_i(t)=e^{B_i(t) - f(t)}, i=1,2$. Max-stable processes that are constructed from log-normal Gaussian spectral processes are commonly referred to in the literature as Brown-Resnick max-stable processes.\\
   By \eqref{foll} if $B_i,i=1,2$ are separable with locally bounded sample paths such that \eqref{aq} holds,
   \kk{then}
   using \eqref{tim} we obtain
\bqn{\label{compSl}
	\pk*{ \sup_{t\in K} X_1(t)> x} \ge  \pk*{ \sup_{t\in K} X_2(t)>x},\quad \forall x>0.
}
Both  processes $Z_1$ and $Z_2$ are Brown-Resnick stationary if \kd{defined} by
$$
Z_i(t)=e^{B_i(t)- \sigma^2_i(t)/2}, \quad i=1,2, \quad t\in \R,
$$
where $B_i,i=1,2$ are two centered Gaussian processes with stationary increments and variance functions $\sigma_1^2$ and $\sigma_2^2$, respectively. Since in general  $\sigma_1$ is different from $\sigma_2$ we cannot use the refinement of Vitale \cite{vitale2} to Slepian inequality stated in \eqref{tim} to arrive at \eqref{compSl}.

\def\WW{W}
\def\xiW{\xi_W}
\def\kalT{E}
 \def\FRE{\mbox{Fr\'{e}chet }}
\def\YD#1{Y^\delta{(#1)}}
\def\xiWD#1{\xi_W^\delta ( #1 ) }
\def\XD#1{X^\delta ( #1 ) }
\def\WD#1{W^\delta ( #1 ) }
\def\XWD#1{X^\delta ( #1 ) }
\def\HWD{\mathcal{H}_W^\delta}
\def\HWDA{\mathcal{H}_W^0}
\def\HWDAA{\mathcal{H}_W}
\def\me{M^\eta}
\def\md{M^\delta}

\def\EID{\theta_W^\delta}
\def\EIDC{\widetilde{\theta_W^\delta}}
\def\EIDD{\widehat{\theta_W^\delta}}
\def\EIDA{\theta_W}

Our next result states the Slepian inequality for Brown-Resnick \bE{max-stable} processes $X_1$ and $X_2$ \bE{which are stationary}. Moreover, it implies a comparison \bE{criterium}  for the corresponding Pickands constants denoted by $\mathcal{H}_{1}$ and $\mathcal{H}_{2}$, respectively.
\kk{Since} the law   of $X_i$'s depends only on their variograms
$\gamma_i(t)=Var(B_i(t)-B_i(0))$, $i=1,2$, \kk{then we} suppose without loss of generality that $\sigma_i(0)=0$ for $i=1,2$.

\BT\label{th.pick}
Let $X_1,X_2$ be two stationary  max-stable Brown-Resnick processes with spectral processes $Z_1$ and $Z_2$, respectively.  Suppose that $\sigma_1(0)=\sigma_2(0)=0$ and $Z_i,i=1,2$ are separable with locally bounded sample paths. If further for any $t\inr$
\bqn{\label{xhg}
	\sigma_1(t)&\ge& \sigma_2(t),% \quad  \forall t\in \R.
}
then for any compact set $K\subset \R$
\bqn{	\pk*{ \sup_{t\in K} X_1(t)> x} \ge  \pk*{ \sup_{t\in K} X_2(t)>x},\quad \forall x>0
	\label{zkhr}
}
and $\mathcal{H}_{1}\ge\mathcal{H}_{2}.$
\ET
\COM{
	\BRM i) \kd{As noted in \cite{WangStoev},} it is not known if there exists
\kd{any} centered Gaussian process $B$  with stationary increments and variance function $\sigma^2$ such that
$\pk{\intT e^{ B(t)- \sigma^2(t)/2} \lambda(dt) < \IF}\in (0,1)$. \\
ii) Under the setup of \netheo{th.pick}, if
$\pk{\mathcal{S}(Z_2)< \IF}>0$, then by Example 1, $\mathcal{H}_{Z_2}>0$ and consequently $\mathcal{H}_{Z_1}>0$. This implies that $X_1$ cannot be generated by a non-singular conservative flow, therefore we can conclude that
$$\pk*{\mathcal{S}(Z_2) < \IF}>0 \Rightarrow \pk*{\mathcal{S}(Z_1)< \IF}>0.$$
iii) \netheo{th.pick} extends the findings of \cite{debicki2002ruin}[Thm 3.2].\\
iv) Slepian inequality has been extended to non-Gaussian setup in \cite{GennaSlepian}.
Using the results of the aforementioned paper, Slepian inequality for max-stable processes can also be derived when $Z_1$  and $Z_2$ have  fidi's satisfying the conditions of \cite{GennaSlepian}[Thm 3.1].

\ERM

\subsection{Wills functional of Gaussian processes and asymptotic constants}
In this section, motivated by Vitale \cite{Vitale1}, we discuss Wills functional for Gaussian processes
and its relation with  Pickands and Piterbarg constants.\\
Consider $Z(t)= e^{B(t)- \sigma^2(t)/2}, t\in \R$, where $B$ is a centered
separable Gaussian process with variance function $\sigma^2$ and bounded sample paths. In \cite{Vitale1} the Wills functional
is defined as
\begin{eqnarray}
 \mathcal{W}_Z(T)= \E*{ \sup_{t \in   [0,T]} {Z(t)}}.\label{Wills1}
 \end{eqnarray}
This functional is important for derivation of lower bounds on supremum of Gaussian processes and random fields. Specifically, by
\cite{Vitale1}[Thm 1] for any $T>0$ we have
\begin{eqnarray}
  \ln \mathcal{W}_Z(T) \le \E[\Big]{\sup_{t \in   [0,T]} B(t) } < \IF,
   \label{wills.bound}
\end{eqnarray}
where the last inequality follows by the boundedness of sample paths, see e.g., \cite{MR1088478}[(3.1)].
Our findings on Pickands constants, when  $Z$ is Brown-Resnick stationary
can be utilised to derive asymptotic
\kd{properties of Wills functionals $\mathcal{W}_Z$
as well as}
lower bounds for $\E{\sup_{t \in   [0,T]} B(t) }$. Indeed, the  Pickands constant $\thZ$ is given in terms of Wills functional as
\begin{eqnarray}
\limit{T} \frac{\mathcal{W}_Z(T)}{T}= \thZ \in (0,\IF).\label{Pick}
\end{eqnarray}
As shown in \cite{K2010} $Z$ is Brown-Resnick stationary if and only if $B$ has stationary increments. In view of Example 1, we have:
%arrive at the following result:

\BEL If $B$ is a centered separable   Gaussian process with stationary increments, locally bounded sample paths and variance function $\sigma^2$, then
\bqn{ \label{betaK}
	\liminf_{T\to \IF} \Bigl[\E[\Big]{\sup_{t \in   [0,T]} B(t) }- \ln T \Bigr]
	& \ge&
	%\limit{T} \bigl[\ln \mathcal{W}_Z(T) - \ln T \bigr] =
	\kd{\ln \thZ,}%>0,
}
provided that \eqref{caj} holds.
\EEL

\BRM
%i) If we do not assume \eqref{caj} above, the lower bound in \eqref{betaK} still holds, but it is equal to 0 if $\pk{\JJ =\IF}=1$. \\
 Sharp bounds on the expectation of maximum of Gaussian processes on finite intervals
 are obtained recently in \cite{MR3574693}. Our result above is of interest when
 considering the growth of supremum of $B$ on large increasing intervals $[0,T]$ and $\sigma^2(t)$ is proportional
 to $a \ln t$ for some $a>8$ (recall \eqref{caj} which guarantees the positivity of $\thZ$).
\ERM

\kd{
Wills functionals are tightly related with Borell-type bounds for tail distribution
of suprema of Gaussian processes. We derive below a variant of Borell inequality accommodated to the setup of this section, see
also \cite{Vitale1}[Cor 1].
For a given variance function $\sigma^2$, let in the following
$\sigma^2_T:=\sup_{t\in [0,T]}\sigma^2(t)$.
\BS\label{Borell1}
Suppose that $B$ is a centered separable Gaussian process with bounded sample paths and variance $\sigma^2$.
For any $T>0$ such that $\sigma_T>0$  and
${V}_u(t)=
e^{{B}_u(t)- Var({B}_u(t))/2}$ with
${B}_u(t)=\frac{u}{\sigma^2_T}B(t)$ we have
\bqn{\label{eqBo}
\pk*{\sup_{t\in [0,T]}B(t)>u}\le \mathcal{W}_{{V}_u}
(T) \exp\left(-\frac{u^2}{2\sigma^2_T}\right)
}
for any $u>0$.
\ES
Next, inserting in \eqref{eqBo} $u=\sigma_T^2$ we get a lower bound for Wills functional (\ref{Wills1}), namely
\[
\IF > \mathcal{W}_Z(T)\ge \pk*{\sup_{t\in [0,T]}B(t)>\sigma_T^2}\exp\left(  \frac{\sigma^2_T}{2} \right).
\]
In the special case of $B$ having stationary increments, utilizing properties of Pickands constants,
we arrive at the following corollary to Proposition \ref{Borell1}.
Set below  $Z_\alpha(t)=\exp\left( B_\alpha(t)-\frac{t^\alpha}{2}   \right)$, $\alpha\in (0,2]$,
with $B_\alpha$  a standard fractional Brownian motion with Hurst index $\alpha/2$.
\BK\label{cor1}
Suppose that $B$ is a centered Gaussian process with stationary increments, \Ecc{bounded sample paths} and variance function $\sigma^2$.
If $\sigma^2(t)\le C t^\alpha$ for $t\in [0,T]$, with $C>0,\alpha\in(0,2]$ and $\sigma_T>0$, then for each $u> 0$
\[
\pk*{\sup_{t\in [0,T]}B(t)>u}\le
\mathcal{W}_{Z_\alpha}\left(\frac{C^{1/\alpha} u^{2/\alpha}T }{\sigma_T^{4/\alpha}} \right)
\exp\left(-\frac{u^2}{2\sigma^2_T}\right).
\]
\EK
Note that $\mathcal{W}_{Z_\alpha}(S)$ is finite and asymptotically linear  in $S$ as $S\to\infty$, by (\ref{Pick}).
}

}

\subsection{Piterbarg constants}
For a given Brown-Resnick stationary process $Z$ as in Section 2 we can define the so-called Piterbarg constants by
$$ 	\mathcal{P}_{f}(K)= \E[\Big]{ \sup_{t \in K }e^{ B(t)- \sigma^2(t)/2  - f(t)} }  $$
for a given positive measurable function $f$, $K=[0,\IF)$ or $K=\R$ and $B$ a centered Gaussian processes with stationary increments as in the previous section.  In the literature, Piterbarg constants
 appear naturally in the tail asymptotics of supremum of non-stationary Gaussian processes, see e.g. \cite{Pit96}.\\
   Clearly, the main question that arises is if
	$\mathcal{P}_{f}(K)$ is finite. So far for $f(t)= a\sigma^2(t),t\in K$ the  finiteness of Piterbarg constants is shown using ideas from the double-sum technique of Piterbarg, see e.g.,  \cite{Pit96}. Utilising the properties of max-stable processes, we are able to show the finiteness of
Piterbarg constants for \kk{more general class of functions} $f$. 	
\kdd{In the following proposition  we consider only the case $K=[0,\IF)$; scenario $K=\R$
follows by analogous line of reasoning. We denote next by $\mathcal{W}$ the Will's functional (see e.g., \cite{Vitale1}) defined by
$$\mathcal{W}(s)=\E*{ \sup_{t \in   [0,s]} e^{B(t)- \sigma^2(t)/2}}, \quad s>0.
$$
}
\BT\label{pit} If $B$ is  a centered Gaussian process with stationary increments, {bounded sample paths}
 and variance function $\sigma^2$, then for any {locally bounded} measurable  $f:[0,\IF] \to \R$ such that $f(t)> a \ln t, a>1, t>0$
%\kd{and $\inf_{t\ge0} f(t)>-\infty$,} \footnote{maybe remove this condition and simply assume $f$ is locally bounded?}
we have that $\mathcal{P}_{f}([0,\IF))< \IF$ and moreover
\bqn{ \label{koloso}
	\mathcal{P}_{f} [(0,\IF))\le \inf_{\delta>0} \mathcal{W}(\delta)  \sum_{i=0}^\IF \sup_{i\in [0,1]} e^{- f(i \delta)}< \IF . %\le \mathcal{H}_Z \Bigl(f(0)+ \int_0^\IF e^{- f(x)}\, dx\Bigr),
}
\ET

\BRM
With  $B$ specified in  \netheo{pit}, for any $T>0$
\cite{Vitale1}[Thm 1]  implies
\begin{eqnarray}
\ln \mathcal{W}(T) \le \E[\Big]{\sup_{t \in   [0,T]} B(t) } < \IF.
\label{wills.bound}
\end{eqnarray}
Consequently, by  Example 1
\bqn{ \label{betaK}
	\liminf_{T\to \IF} \Bigl[\E[\Big]{\sup_{t \in   [0,T]} B(t) }- \ln T \Bigr]
	& \ge&
	%\limit{T} \bigl[\ln \mathcal{W}_Z(T) - \ln T \bigr] =
	\kd{\ln \thZ> {- \IF,}}%>0,
}
provided that \eqref{caj} holds.
\ERM

 {
 	\subsection{Growth of supremum}
 	Let $X$ be a separable max-stable stationary process with locally bounded sample paths and
 spectral process $Z$ such that \eqref{eq1} holds.
 \kk{For any $T>0$, let}  $A_T:= M_X(T)T^{-1/\alpha},T>0$.
 \kk{Since} $A_T$ is non-negative, then   for any $\beta\in (0,\alpha)$
 	we have
 \begin{eqnarray}
 	 \E{ A_T^\beta } &=& \int_0^\IF \pk{A_T > x^{1/\beta} } dx\nonumber\\
 		&=& 	\int_0^\IF \Bigl (1- e^{-  \E{ \sup_{t\in [0,T]} Z^\alpha (t) }/(T x^{\alpha/\beta})} \Bigr)dx\nonumber\\
 		&=& \Bigl(\E[\Big]{ \sup_{t\in [0,T]} Z^\alpha (t) }/T\Bigr)^{\beta/\alpha} \Gamma(1- \beta/\alpha),\label{f.13}
 \end{eqnarray}
 	where $\Gamma(\cdot)$ stands for the Euler Gamma function and
 \kk{in (\ref{f.13})} we used \eqref{fred}. Consequently, as $T \to \IF$
 	\bqn{   \E{ A_T^\beta }  \to \thZ ^{\beta/\alpha} \Gamma(1- \beta/\alpha)
 		< \IF \label{limea}.
 	}
% holds as $T\to \IF$.
 	A direct implication of \eqref{at1}, \eqref{at2} and \eqref{limea} is the following result.
 	\BS \label{ProK} If $X(t),t\in \mathcal{T}$ where $\mathcal{T}= \R$ or $\Z$ is  a max-stable stationary process as above, then \eqref{limea} holds and moreover $A_T^\beta,T>0$ is uniformly integrable for any $\beta \in (0,\alpha)$.
 	\ES
 }
 \BRM For $X$ a $S\alpha S$ stationary random field \eqref{limea}  has been shown in \cite{RoyY}[Thm 3.1]. Extension of \eqref{limea} to stationary max-stable random fields is straightforward and omitted here.
 \ERM
 \section{Proofs}

\prooftheo{proviant} We adapt  the  arguments of the proof of \cite{Genna04c}[Thm 2.1] for our max-stable process.
 Note first that by \eqref{bask}, in view of \cite{dom2016} $X$ has a dissipative Rosi\'nski  representation with some
 process $L$, which by the construction in the aforementioned paper is stochastically continuous and locally bounded
 (these properties are inherited from $X$). Therefore \kk{there exists} jointly measurable and separable version of
 $L$ which is locally bounded; we shall consider this version below. \\
\underline{Step 1}: Since $Z$ is given by \eqref{piratenQ}, and moreover $Z$ is locally bounded, then by \eqref{spilman},
for any $T>0$, we have that $\thZ \le \E{M_{Z^\alpha }(T)}/T $.
\kk{Moreover,} by stationary of $X$ %\bE{and \eqref{nat}}
$$ H(T): = \E*{ \sup_{ 0 \le t \le T} Z^\alpha (t)}=
\E*{ \sup_{ 0 \le t \le T} Z^\alpha (-t)} \in \bE{(0, \IF)}.$$
    Consequently, as in \cite{Genna04c} we obtain the following lower bound
    %for $\NN$ with density $p(t)>0, t\in \TT$ being independent of $L$
\bqn{
\IF > H(2) &=&   \E*{ \sup_{0 \le t \le 2} Z^\alpha(-t) }\notag\\
&=&  c\E*{ \sup_{0 \le t \le 2}  \frac{B^\NN L^\alpha(-t) }{p(\NN)} }\notag \\
&=& c \E*{ \intT  \sup_{0 \le t \le 2}   L^\alpha(x-t)  \lambda(dx)} \label{xhde}\\
&=&c\sum_{i\in \Z} \E*{ \int_{i}^{i+1} \sup_{0 \le t \le 2}   L^\alpha(x-t)
\lambda(dx)} \notag \\
&\ge & c\sum_{i\in \Z} \E*{  \sup_{i-1 \le t \le i}   L^\alpha(t)}\int_{i}^{i+1}
\lambda(dx) \notag \\
&\ge &c \E*{  \sup_{t \in \R }   L^\alpha(t)}. \notag
}
Further,
since we assume that
$X(0)$ has \FRE  $\Phi_\alpha$ distribution, by \eqref{foll}
\bqny{%\label{constc}
	\quad \quad - \ln \pk{X(0)\le 1}= \E{Z^\alpha(0)} = c \E*{ \intT (B^t L^\alpha)(0) \lambda(dt)} = 	 c \E*{ \intT L^\alpha(t)\lambda(dt)}>0.
}
\kk{Hence we conclude that }
\bqn{ \label{cabina}
\E*{  \sup_{t \in \R } L^\alpha(t)} \in (0,\IF).
}
\underline{Step 2}:  If for some positive $M$ we have $\pk{ \sup_{\abs{t}\ge M} L(t) =0}=1$, then
 for $T> 2M$ %and $\NN$ as above
\bqny{
	\frac{H(T)}{T}	
	&=& \frac{c}{T} \intT \E*{ \sup_{0 \le s \le T}    L^\alpha (x+s) } \lambda(dx)\\
	&=& \frac{c}{T} \intT \E*{ \sup_{t \le v \le T+t}    L^\alpha (v) } \lambda(dt)\\
	&=& \frac{c}{T} \int_{-M-T}^{M} \E*{ \sup_{t \le v \le T+t}    L^\alpha (v) } \lambda(dt)\\
	&=& c \int_{-M/T}^{1+M/T} \E*{ \sup_{-Tx \le v \le T(1-x)}   \kd{L^\alpha(v) }} \lambda_T(d x)/T\\
	&=&  c(1- 2M/T) \E*{ \sup_{-M \le v \le M}    L^\alpha (v) }
	\int_{0}^{1}  \lambda_T(d x)/T \\
	&&  + O(M/T)\\
	&\to& c\E*{ \sup_{v\inr }    L^\alpha (v) }, \quad T\to \IF,
}
where $\lambda_T(dx)= \lambda(T dx)$.\\
\underline{Step 3}:  For $M>0$ set $L_M(t)= L(t) \mathbb{I}(\abs{t}\le M)$. Using \eqref{cabina} and applying the bounded convergence theorem we obtain
\bqn{\label{aufw}
	\limit{M} \Delta_M:= \limit{M}  \Bigl[\E*{ \sup_{t \in   \R} L_M^\alpha(t)} - \E*{ \sup_{t \in   \R} L^\alpha(t)}\Bigr]=0.
}	
Further by \eqref{xhde}, again the bounded convergence theorem yields  %and %for any $T>0$
\bqn{\label{zug}
	\limit{M} \E*{ \intT \sup_{t \in   [0,1]} [L^\alpha(t-x) -L_M^\alpha(t-x)]   \lambda(dx) }=0.
}
Write next for $T>0$
\bqny{
\kk{\Abs{ \frac{H(T)}{T}- \thZ}}&\kk{\le}&
\frac{1}{T}\Abs{ H(T)-c\E*{ \intT  \sup_{0 \le s \le T}    L_M^\alpha (x+t)
	 \lambda(dx) } }\\
&& +
c\Abs{ \frac{1}{T} \E*{ \intT  \sup_{0 \le t \le T}    L_M^\alpha(x+t)
	\lambda(dx)} - \E*{ \sup_{t \in  \R} L_{M} ^\alpha (t)}} + c\abs{\Delta_M}.
}
The second and third terms converge as $T\to \IF$ and $M\to \IF$, respectively by Step 2 and \eqref{aufw}.
Further by \eqref{zug} (write $[T]$ for the smallest integer larger than $T$)
\bqny{
\lefteqn{\frac{1}{T}\Abs{ H(T)-c\E*{ \intT  \sup_{0 \le t \le T}    L_M^\alpha(t-x) \lambda(dx)}} }\\
	&\le &
	\frac{c}{T}\E*{ \intT  \sup_{0 \le t \le T} [ L^\alpha (t-x) -
		L^\alpha_M (t-x) ] \lambda(dx) }\\
	&\le &
\frac{c}{T}\E*{ \intT  \sum_{ 1 \le j \le [T]} \sup_{j-1\le t \le j} [ L^\alpha (t-x) -
	L^\alpha_M (t-x) ] \lambda(dx) }\\
	&\le &
\frac{c[T]}{T}\E*{ \intT  \sup_{0 \le t \le 1} [ L^\alpha (t-x) -
	L^\alpha_M (t-x) ] \lambda(dx) }\\
&\to& 0, \quad M\to \IF,
}
hence the proof is complete.
\QED

\prooftheo{th2} We consider first the case that $\pk{\JJ< \IF}=1$. Since $Z$ has c\`adl\`ag sample paths we can assume without loss of generality that $Z$ is such that  $\JJ>0$ almost surely. 
%Note that  $\E{Z^\alpha(t)}=1, t\inr$ implies that
%$\JJ >0$ almost surely, since by Fubini theorem
%$$\E{\JJ} \ge \E*{\int_a^b Z^\alpha(t)\lambda(dt) }=
%\int_a^b \E{ Z^\alpha(t)}\lambda(dt)= b- a.$$
Hence %\bE{in view of \eqref{nat}}
almost surely $\JJ\in (0,\IF) $.\\ %and $Z(0)\in [0, \IF)$.\\
The map $H: D \to [0,\IF]$ defined by $H(f) \to \intT \abs{f(t)}^\alpha\lambda(dt)=:\mathcal{S}(f)$ is $\mathcal{D}/\mathcal{B}(\R)$ measurable since
$\mathcal{S}(f)$ for $f$ c\`adl\`ag is determined by $f(t), t\in \mathcal{Q}$ with $\mathcal{Q}$ the set of all rational \bE{numbers}, and  $\mathcal{D}$ coincides with the $\sigma$-algebra  $\sigma(\pi_t,t\in \mathcal{Q})$, see \cite{MR838085}[Prop 7.1].  Hence we have 	
\bqn{
	1&=&\E{ Z^\alpha(0)}= \E{Z^\alpha(0) \mathbb{I}( \JJ \in (0,\IF))} \notag\\
	&=& \E{ Z^\alpha(0) \mathbb{I}(0<  Z(0)< \IF,  \JJ \in (0,\IF))} \notag\\
	&=& \E{ Z^\alpha (0)\mathbb{I}(\mathcal{S}(Z/Z(0)) \in (0,\IF))} \notag\\
	&=& \E{ \mathbb{I}(\mathcal{S}(\Theta) \in (0,\IF))}
	\label{leeds}
}
\COM{and and similarly,
\bqny{
	1&=&\E{ Z^\alpha(0)}= \E{Z^\alpha(0) \mathbb{I}( \limit{t} Z(t)=0 )} \\
&=& \E{\limit{t} \Theta(t)=0 )}
}
thus since by \cite{dom2016} $\pk{ \JJ\in (0,\IF)}=1$ if and only if $\pk{\limit{t} Z(t)=0}=1$, we conclude   that
}
implying
\bqn{ \pk{ \JJ\in (0,\IF)}=1 \Leftrightarrow \pk{\mathcal{S}(\Theta) \in (0,\IF)}=1.
}	
Consequently,  the random shape function $L$ given by
$$L(t)= \frac{\Theta (t)}{(\mathcal{S}(\Theta))^{1/\alpha}},\quad t\inr$$
  is well-defined with   c\`adl\`ag sample paths and $\mathcal{S}(L)=1$ almost surely.\\
 We continue by showing that   $\tZ:= (p(\NN))^{-1/\alpha}B^\NN L$ is a spectral process such that its corresponding max-stable process $\tX$ with de Haan representation   \eqref{eq1} (taking $\tZ$ instead of $Z$) is stationary.  Since for any $h\in \R$
 $$ B^h \tZ = (p(\NN))^{-1/\alpha}B^{\NN+h} L =
 (p_h(\NN_h))^{-1/\alpha}B^{\NN_h} L =:Z_{N_h},$$
 with $\NN_h=\NN +h$ which has density function $p_h(t)=p(t-h)$ it follows that $\tX$ is stationary, since $Z_{N_h}$ is a spectral process for $\tX$ by the shift-invariance of Lebesgues measure. Next, we prove that
$\tX$ has the same fidi's as $X$. By the stationarity, in view of \eqref{argmax2} this follows if we show that the spectral tail process $\tTh$ of $\tX$ has the same fidi's as $\Theta$.
%Finally, we show that $\tX$ corresponding to $\tZ = (p(\NN))^{-1/\alpha} B^\NN \tL $ has fidi's equal to $X$.
  For any $A \in \mathcal{A}$ %(set $U=   B^\NN \Theta/\Theta(0))$
\bqny{ \pk{\tTh \in A}	&=& \E*{ \tZ^\alpha(0) \mathbb{I}( \tZ/ \tZ(0) \in A)}
\\
&=& \E*{   \frac{ ( B^\NN \Theta) ^\alpha(0)  }{\mathcal{S} (\Theta)}\mathbb{I}(  (B^\NN \Theta)/ (B^\NN \Theta)(0)\in A, \JJ \in (0,\IF) )}\\
&=&  \intT \E*{  \frac{ (B^t \Theta)^\alpha(0)  }{\mathcal{S}  (B^t \Theta)}\mathbb{I}( (B^t \Theta)/(B^t \Theta)(0)\in A, \mathcal{S}(B^t\Theta) \in (0,\IF))} \lambda(dt)\\
&=&  \intT \E*{  H(B^t\Theta)} \lambda(dt).
}
The functional $H(f)= F(f)$ is $0$-homogeneous non-negative and
$\mathcal{D}/\mathcal{B}(\R)$ measurable
(\kd{we use the convention}  $0\IF=:0$).
Since $\mathcal{S}(B^t f)= \mathcal{S}(f),t\inr$, \kd{then} applying \nelem{sebi} in Appendix we obtain
\bqny{ \intT \E*{  H(B^t\Theta)} \lambda(dt) &=&  \intT \E*{ Z^\alpha (t) H(Z)} \lambda(dt)\\
	&=& 	 \intT \E*{  Z^\alpha(t)\frac{ Z^\alpha(0)}{\JJd} \mathbb{I}( Z/Z(0)\in A, \JJ \in (0,\IF)} \lambda(dt)\\
&=& 	 \E*{ \frac{ \intT   Z^\alpha(t)\lambda(dt) }{\JJd}  Z^\alpha(0) \mathbb{I}( Z/Z(0)\in A, \JJ \in (0,\IF)} \\
&=&  \E*{ Z^\alpha(0)\mathbb{I}( Z/Z(0)\in A)} \\
&=&  \pk{\Theta \in A},
}
establishing \bE{that $ X$ has the same fidi's as $\tX$} and the dissipative Rosi\'nski  representation of $X$ with $L$ constructed above. \\
 {\bE{Next}, for a given spectral process $Y$ we denote by $\mathcal{H}_Y$ the corresponding Pickands constant.}
 Next,  if $\pk{\JJ< \IF}\in (0,1)$ by \eqref{zetC} we have
 $$ \mathcal{H}= a \limit{T} \E*{ \sup_{t \in   [0,T]} Z_C^\alpha (t)/a  } = a \mathcal{H}_{Y},$$
 where $Y= Z_C/a^{1/\alpha}$ with
 $Z_C(t)=
Z(t) \mathbb{I}( 0<\JJ< \IF)$ and $a= \pk{0< \mathcal{S}(\Theta) <  \IF }\in (0,1)$. Note in passing that
$$ \E{Y^\alpha(0)} =\E{Z_C^\alpha(0)/a}=
\E{Z^\alpha(0) \mathbb{I}(0<\JJ< \IF)}/a =1.$$
 The spectral tail process $\Theta_C$ of $Y $ is calculated for any $A \in \mathcal{A}$ by
\bqny{ \pk{\Theta_C \in A } &=&
	 \E*{ \frac{Z^\alpha_C(0)}a \mathbb{I}( Z_C/ Z_C(0) \in A) }\\
	 &=& 	 \E*{ \frac{Z^\alpha(0)}a \mathbb{I}( Z/ Z(0) \in A) \mathbb{I}(0<\JJ< \IF ) }\\
	 &=& 	 \E*{ \frac{ \mathbb{I}( \Theta  \in A)}{a} \mathbb{I}(0<\JJ< \IF ) }\\
 &=& \pk{\Theta \in A \lvert 0 < \mathcal{S}(\Theta)< \IF}.
}
As above for the stationary max-stable  process $X_C^*$ with spectral process $Y $ we have that it has a dissipative Rosi\'nski representation with random shape function $L_C$ given by
$$L_C(t)= \frac{\Theta_C(t)}{\mathcal{S} (\Theta_C)}
= \frac{\Theta(t)}{\mathcal{S}  (\Theta) } \Bigl \lvert (0 < \mathcal{S}(\Theta)< \IF),\quad t\inr .$$
Consequently,
%since  \eqref{zetC} we have
%$\thZ= \mathcal{H}_{Z_C}$ then
using further \netheo{proviant}
\bqny{
	 \thZ  	&=& a  \mathcal{H}_{Y}\\
	  &=&
	a  \frac{ \E{ \sup_{t \in   \R}  L^\alpha_C (t) } }
	 { \E{\mathcal{S}(L_C) }}\\
&=& 	a \frac{ \E*{ \sup_{t \in   \R}  \Theta_C^\alpha (t)/\mathcal{S}(\Theta_C) } }
{ \E{\mathcal{S}( \Theta_C/ \mathcal{S}(\Theta_C)) }}\\
&=& 	
% \E*{ \sup_{t \in   \R}  \frac{\Theta_C^\alpha (t)}{ \delta \mathcal{S}_\delta(\Theta_C)} } =
a \E*{ \sup_{t \in   \R}  \frac{\Theta_C^\alpha (t)}{  \mathcal{S}(\Theta_C)} } \\
&=& \E*{ \sup_{t \in   \R}  \frac{\Theta^\alpha (t)}{  \mathcal{S}(\Theta)};0< \mathcal{S}(\Theta)< \IF }
}
establishing the proof.    \QED

 \def\Cov{\mathrm{Cov}}

%\proofprop{pr2}
 \prooftheo{th.pick} In view of \cite{KabExt} both $X_1$ and $X_2$ are max-stable statoinary processes. Hence in view of \eqref{foll},
 in order to show \eqref{zkhr}  we need to prove  that for any compact set $K \subset \R$
\begin{eqnarray}
\E*{\sup_{t\in K}Z_2(t)} \le \E*{ \sup_{t\in K} Z_1(t)}<\infty
 \label{Pick_b}
\end{eqnarray}
\bE{is valid} with $Z_i(t)= e^{B_{i}(t)- \sigma_i^2(t)/2}$, $i=1,2$.  \bE{Note in passing that the finiteness} of  $\E*{ \sup_{t\in K} Z_1(t)}$ follows from \cite{Vitale1}[Thm 1]. \\
 By the  stationarity of increments of $B_i$'s the variance functions  $\sigma^2_i(t),t\inr, i=1,2$ are negative definite functions.
Consequently, by Schoenberg theorem, for each $u>0$, $i=1,2$
\kd{function}
 \[R_u^{(i)} (s,t):=\exp \left(- \frac{1}{2 u^2}\sigma^2_i(s-t)\right), \quad s,t\inr
 \]
 is positive definite and thus a valid covariance function.
 %\footnote{do we need the continuity of sample paths, maybe we assume only separability?}

 Let   $W_u^{(i)}(t), t\in\R$, $u>0$ be a family of separable
 centered stationary Gaussian processes with covariance functions
 \[\Cov(W_u^{(i)}(s),W_u^{(i)}(t))=R_u^{(i)} (s,t), \quad s,t\inr \]
 for $i=1,2$. Since by assumption $\sigma_1(t) \ge \sigma_2(t)$ for any $t\in \R$,
 then for any  $s,t\in  \R$ we have
 $$  R_u^{(1)} (s,t) \le R_u^{(2)} (s,t).$$
 Hence, for a given compact set $K\subset \R$,
 applying Slepian inequality, see, e.g., \cite{MR1472736}[Thm 3] %  and \cite{GennaSlepian} for an extension to non-Gaussian case),
for any $u>0$
 \begin{eqnarray}
 \pk*{\sup_{t\in K}W^{(1)}_u(t)>u}\ge \pk*{\sup_{t\in K}W^{(2)}_u(t)>u}.\label{slepian}
 \end{eqnarray}
The definition of the covariance functions above yields  for $i=1,2$
 \[
 \lim_{u\to\infty} \sup_{s,t\in K} \left|         \frac{1-\Cov(W_u^{(i)}(s),W_u^{(i)}(t))}{\frac{1}{2 u^2}\sigma^2_i(s-t)} -1 \right|=0.
 \]
 Applying \cite{DeHaJi2016}[Lem 6.1] %\cite{DHL16}[Theorem 2.1]
 for any $\eta>0$ and any given compact $K \subset \R$ %$a<b$ constants %K=\prod_{i=1}^d [a_i,b_i],$  $a_i< b_i$
 \bE{we obtain}
%\begin{eqnarray}
% \lim_{u\to\infty}\frac{\pk{\sup_{t\in \eta \Z^d\cap P}\frac{W^{(i)}_u(t)}{1-f(t)/u^2}>u}}{\Psi(u)}=
% \E*{ \sup_{t\in \eta \Z^d\cap P}{e^{B_{i}(t)- \sigma_i^2(t)/2 +f(t)}}},\label{from.uniform}
% \end{eqnarray}
 \begin{eqnarray}
 \lim_{u\to\infty}\frac{\pk{\sup_{t\in \eta \Z\cap K }W^{(i)}_u(t)>u}}{\Psi(u)}=
 \E*{ \sup_{t\in \eta \Z\cap K } Z_i(t)}.\label{from.uniform}
 \end{eqnarray}
% where $Z_i(t)= \kd{e^{B_{i}(t)}- \sigma_i^2(t)/2}$.
 Hence, by the separability and local boundedness of the sample paths, (\ref{slepian})
 combined with (\ref{from.uniform}) implies \eqref{Pick_b}.\\
The assumption
$\E{Z_i(t)}=1$ for any $t\inr$ implies that
  the max-stable processes corresponding to $Z_1$ and $Z_2$  have unit
  \FRE marginals $\Phi_1$.   The Pickands constants corresponding to $Z_1$ and $Z_2$  denoted by $\mathcal{H}_1$ and $\mathcal{H}_2$, respectively exist and are finite since $X_1$ and $X_2$ are stationary (see the argument given for the derivation of \eqref{spilman}). Hence a direct application of \eqref{Pick_b} implies that $\mathcal{H}_1 \ge \mathcal{H}_2$. This completes the proof.
\QED

\Ecc{\prooftheo{pit} First note that by the fact that $Z(t)=e^{B(t)- \sigma^2(t)/2}$ is Brown-Resnick stationary, we have that
	$$ \E[\Big]{ \sup_{t\in [\delta i, \delta (i+1)]} Z(t)}=% e^{B(t)- \sigma^2(t)/2}}=
	\E[\Big]{ \sup_{t\in [0, \delta]}Z(t)}=:% e^{B(t)- \sigma^2(t)/2}}=:
	\mathcal{W}
	(\delta) \bE{\in (0,\IF)}$$
	for any \bE{$\delta \inr,i>0 $}.
	Hence for any  positive $\delta$ using the assumption that $f$ {is locally bounded} and $f(t)> a \ln t$ for all $t$ large with some $a>0$,  we obtain
	\bqny{
		\E[\Big]{ \sup_{t\ge 0 }Z(t) e^{-f(t)}} %e^{B(t)- \sigma^2(t)/2 - f(t)}}
			& \le &
		\sum_{i=0}^{\IF }\E[\Big]{ \sup_{t\in [\delta i, \delta (i+1)]}Z(t)%e^{B(t)- \sigma^2(t)/2
			e^{- f(t)}}\\
&\le & 		\sum_{i=0}^{\IF }\E[\Big]{ \sup_{t\in [\delta i, \delta (i+1)]}Z(t)}%e^{B(t)- \sigma^2(t)/2
\sup_{t\in [\delta i, \delta (i+1)]}	e^{- f(t)}\\
		&= & \mathcal{W} (\delta) \sum_{i=0}^{\IF} \sup_{t\in [i \delta, (i+1) \delta]}e^{f(t)} < \IF,
	}
	hence the proof follows. 	
	\QED
}

\section{Appendix: Tilt-Shift and inf-argmax fomula}
Next we present the tilt-shift formula which is initially shown for the special case of Brown-Resnick
max-stable processes with log-normal $Z$ in \cite{DM}.  The inf-argmax formula mentioned above is shown in \cite{Htilt}, we present below a shorter proof. %Note that tilt-shift formula is equivalent with time-change formula of \cite{BojanS} for the discrete case, see e.g., \cite{Htilt}.
\BEL \label{sebi} Let $X(t),t\inr$ be a max-stable process with \FRE marginals $\Phi_\alpha$, de Haan  representation \eqref{eq1} and c\`adl\`ag sample paths.\\
i)  If $X$ is stationary, then for any non-negative  $0$-homogeneous $\mathcal{D}/\mathcal{B}(\R)$-measurable functional $H$  we have
\BQN \label{stgh}
\E{ Z^\alpha(h) H(Z) } = \E{ H(B^h \Theta)}= \E{Z^\alpha(0)H(B^h Z)}.
\EQN
ii) If \eqref{stgh} holds for any $h\in \TT$, then $X$ with representation \eqref{eq1} is stationary.
\EEL
\prooflem{sebi} $i)$  As in \cite{SBK} we have that the stationarity of $X$ implies the shift invariance of the exponent measure, i.e.,  for any $h\in \TT,A \subset \mathcal{D}$
$$ \nu(A)=\int_0^\IF \pk{ uZ \in A} \alpha u^{-\alpha-1} \,du =
\int_0^\IF \pk{ u B^h Z \in A} \alpha u^{-\alpha-1} \,du.$$
If $A$ is $0$-homogeneous (meaning $cA=A, c>0$) set in $\mathcal{D}$, then
since further $\E{Z^\alpha(h)}=1$ implies that $\pk{Z^\alpha(h) \in [0,\IF)}=1$
for any $h\in \TT$ using $z^\alpha= \int_0^\IF \mathbb{I}( r z >1)  \alpha r^{-\alpha-1} \,dr$ valid for any $z\in (0,\IF)$ we obtain
\bqny{
	\E{Z^\alpha (h) \mathbb{I}(Z \in A)} &=&
		\E*{Z^\alpha (h) \mathbb{I}(Z \in A, Z(h) \in (0,\IF)}  \\
&=&
	\E*{ \int_0^\IF  \mathbb{I}( r Z(h) >1)\mathbb{I}( rZ \in A)  \alpha r^{-\alpha-1} \,dr } \\
	&=&  	\E*{ \int_0^\IF  \mathbb{I}( r Z(0) >1, rB^h Z \in A) \alpha r^{-\alpha-1} \,dr } \\
	&=&  	\E*{ \int_0^\IF  \mathbb{I}( r Z(0) >1)  \alpha r^{-\alpha-1} \,dr \mathbb{I}( B^h Z \in A)} \\
	&=&  	\E*{ Z^\alpha(0)  \mathbb{I}( B^h Z \in A) }\\
	&=&  	\E*{ Z^\alpha(0)  \mathbb{I}( B^h Z/Z(0) \in A) }\\
	& =&
	\E*{   \mathbb{I}( B^h \Theta \in A) },%= 	\pk{    B^h \Theta \in A},
}
hence the claim for  any $0$-homogeneous $\mathcal{D}/\mathcal{B}(\R)$-measurable functional $H$ follows easily.\\
ii)  If \eqref{stgh} holds, then for any $n\ge 1, t_i \inr, x_i>0, i\le n$
since $\inf \arg \max $ functional is $0$-homogeneous and  $\mathcal{D}/\mathcal{B}(\R)$-measurable, by \eqref{foll} we have
\bqny{ - \ln \pk{ X(t_i) \le x_i, 1 \le i\le n}&=&
	\E*{ \max_{1  \le i \le n}  x_i^{-\alpha}  Z^\alpha (t_i)  }\\
%	&=& \sum_{k=1}^n 	\E*{    \max_{1  \le i \le n} Z^\alpha (t_i)/x_i^\alpha ; \inf \arg \max_{1  \le i \le n} Z^\alpha(t_i)/x_i^\alpha =k}\\
	&=& \sum_{k=1}^n x_k^{-\alpha}	\E*{    Z^\alpha (t_k)\mathbb{I} \Bigl(\inf \arg \max_{1  \le i \le n} Z^\alpha(t_i)/x_i^\alpha=k\Bigr)}\\
	&=:& \sum_{k=1}^n x_k^{-\alpha}	\E*{    Z^\alpha (t_k) F_k(Z)}\\
	&=& \sum_{k=1}^n x_k^{-\alpha}	\E*{  F_k(B^{t_k} \Theta)}\\
	&=& \sum_{k=1}^n x_k^{-\alpha}	\pk*{\inf \arg \max_{1  \le i \le n} \Theta^\alpha(t_i-t_k)/x_i^\alpha=k }, }
where we used \eqref{stgh} in the second last line above.
Consequently, $X$ is stationary and thus the proof if complete.
 \QED

	\section*{Acknowledgments}
	Many thanks to  Parthanil Roy for discussions and suggestion of the key reference \cite{Genna04c}. We thank the referees for  numerous suggestions that improved the original manuscript. 	EH was supported by SNSF Grant 200021-175752/1.  KD
was partially supported by NCN Grant No 2015/17/B/ST1/01102 (2016-2019).

\bibliographystyle{ieeetr}
\bibliography{EEEA}
\end{document}